\documentclass[a4paper]{amsart}
\usepackage{amsmath}
\usepackage{amssymb}
\usepackage{bbm}
\usepackage{a4wide}
\usepackage{graphicx}
\usepackage{color}

\parskip\smallskipamount

\def\<#1>{\langle#1\rangle}

\def\e{\mathrm{e}}
\def\e1{\varepsilon}

\begin{document}

 \author[Shalosh B. Ekhad]{Shalosh B. Ekhad}
 \address{Shalosh B. Ekhad c/o Doron Zeilberger, Mathematics Department, Rutgers University (New Brunswick), Piscataway, NJ, USA.}
 \email{c/o zeilberg@math.rutgers.edu}

 \author[Amitai Regev]{Amitai Regev}
 \address{Amitai Regev, Department of Mathematics, Weizmann Instiute of Science, Rehovot, Israel}
 \email{ amitai.regev@weizmann.ac.il}

\date{July 30, 2010}

\title[Refined Asymptotics and Explicit Recurrences for the numbers of Young tableaux in the $(k,l)$ hook for $k+l \leq 5$]
{Refined Asymptotics and Explicit Recurrences for the Numbers of Young tableaux in the $(k,l)$ hook for $k+l \leq 5$}
 \maketitle

\section*{The Input: the Sequences $S_{k,l}^{(z)}(n)$}

Recall that a {\it Young  Shape} (alias {\it Ferrers diagram}) $\lambda=(\lambda_1, \lambda_2, \dots )$
is an arrangement of 
left-justified  rows of boxes, with $\lambda_1$ boxes in the first row,
$\lambda_2$ boxes in the second row, etc. If the total number of boxes is $n$, we write
$\lambda \vdash n$.

The set of shapes $\lambda$ with $n$ cells with the property that $\lambda_{k+1} \leq l$ is denoted
by $H(k,l;n)$.

Recall also that a {\it Standard Young Tableau} (SYT) of shape $\lambda \vdash n$ is any way of placing the
integers from $1$ to $n$ inside the $n$ boxes in such a way that all rows and all columns are increasing.
There is a well-known, explicit, formula for $f^{\lambda}$, the number of SYTs of a given shape $\lambda$,
due to Young and Frobenius, as well as a more elegant reformulation, called the
hook-length formula, due to Frame, Robinson, and Thrall.

In \cite{BR}, Regev and Berele studied
the sequences
$$
S_{k,l}^{(z)}(n)=
\sum_{\lambda \in H(k,l;n)} (f^{\lambda})^{z} \quad ,
$$
and obtained explicit leading asymptotics for {\it general} $(k,l)$
as well as general (real) $z$.

Except for $(k,l) \in \{ (1,0),(2,0),(4,0),(1,1)\}$ (when $z=1$) 
and $(k,l) \in \{ (1,0),(2,0),(1,1), (\infty, \infty) \}$ (for $z=2$),
there are no known ``nice'' expressions for $S_{k,l}^{(z)}(n)$ in terms
of $n$. Of course, using the Young-Frobenius formula one can always
express it as a {\it multi-sum} of hypergeometric terms, but, in general
with $k+l-1$  sigma signs.

 The next best thing to a closed-form formula is a {\it linear recurrence}.
It follows from the Fundamental Theorem of Wilf-Zeilberger theory (\cite{WilZ})
that for any specific $(k,l)$ and {\it positive integer} $z$, the sequence
$A(n)=S_{k,l}^{(z)}(n)$ is {\it holonomic} (alias $P$-recursive), in other words it satisfies
a (homogeneous) linear recurrence equation with {\it polynomial} coefficients. This 
means that there exists a positive integer $L$, and polynomials $p_0(n), \dots, p_L(n)$ such that
$$
\sum_{i=0}^{L} p_i(n) A(n+i) =0 \quad , \quad n \geq  0 \quad .
$$

Furthermore there exist algorithms for computing these recurrences, but they are rather slow for
multi-sums.

Since we know {\it a priori} that such a recurrence exists, we may just as well find it {\it empirically}
by asking the computer to compute the first $60$ (or whatever) terms and then
guess the recurrence by {\it undetermined coefficients}. This is done by the command
{\tt listtorec} in Bruno Salvy and Paul Zimmermann's Maple package {\tt gfun}
(but we prefer to have our own version).

While there is a ``theoretical'' possibility that the guessed recurrence is not the ``right'' one,
this is very unlikely, and besides, using the recurrence one can extend the sequence very fast,
and check whether the continuation matches, and if it does 
(say for the next $20$ terms), then
this is a {\it semi-rigorous} proof, good enough for us. We know that we have the {\it option} to
find a fully rigorous proof by the {\it multivariable Zeilberger algorithm} \cite{AT}
(that alas, will take much longer, and for larger $(k,l)$ may not terminate in a reasonable
amount of time) but {\it why bother?}.

The beauty of knowing a linear recurrence is that one can use the  Poincar\'e-Birkhoff-Trjitzinsky
method to derive higher-order asymptotics for the sequence. This method, described in \cite{WimZ},
is fully implemented in Doron Zeilberger's Maple package

\verb|http://www.math.rutgers.edu/~zeilberg/tokhniot/AsyRec| \quad .

The drawback of that method is that it does not output the {\it constant} factor,
that must be determined empirically, in general. Luckily, 
thanks to \cite{BR}, we know
that constant exactly (expressed in terms of $\pi$), and this enables the computation
of very refined asymptotics of the sequences 
$S_{k,l}^{(z)}(n)$ for the most interesting cases $z=1$ and $z=2$, and for
$k+l \leq 5$.

\section*{The Output: Recurrences and asymptotics for $S_{k,l}^{(1)}(n)$ and $S_{k,l}^{(2)}(n)$
for $k+l \leq 5$}

All the necessary algorithms have been implemented in the self-contained Maple package
{\tt HOOKER}, kindly written by Doron Zeilberger, 
and available from

\verb|http://www.math.rutgers.edu/~zeilberg/tokhniot/HOOKER| \quad .

\bigskip

Detailed output files, with lots of terms in the sequences, and asymptotics to order $10$, can be viewed in the webpage of this article:

\verb|http://www.math.rutgers.edu/~zeilberg/mamarim/mamarimhtml/hooker.html| \quad .

Here we only list (for the record)
the output for $(k,l) \in \{(2,1),(2,2)\}$ and $z=1,2$, and even for these, we only list the
first few terms, the annihilating operator and the refined asymptotics to order 3.

Below $N$ denotes the shift-operator in $n$, so, for example, the annihilating operator for
the Fibonacci sequence is $N^2-N-1$ and for $n!$ is $N-n-1$.

{\bf $(k,l)=(2,1)$, $z=1$}

$$
1, 2, 4, 10, 26, 71, 197, 554, 1570, 4477, 12827, 36895, 106471, 308114, 893804, 2598314, 7567466, 22076405, 
$$
$$
64498427, 188689685
$$
$$
3\,{\frac {n+2}{n+3}}-{\frac {n \left( n+2 \right) N}{ \left( n+3 \right)  \left( n+1 \right) }}-{\frac { \left( 9+11\,n+3\,{n}^{2} \right) {N}^{2}}{
 \left( n+3 \right)  \left( n+1 \right) }}+{N}^{3}
$$

$$
1/4\,{3}^{n}\sqrt {{n}^{-1}} \left( 1-3/16\,{n}^{-1}+{\frac {1}{512}}\,{n}^{-2}+{\frac {135}{8192}}\,{n}^{-3} \right) \sqrt {3}{\frac {1}{\sqrt {\pi }}}
$$

{\bf $(k,l)=(2,1)$, $z=2$}

$$
1, 2, 6, 24, 120, 695, 4403, 29540, 206244, 1483371, 10919271, 81896661, 623810421, 
$$
$$
4813777566, 37561178658, 295907998908, 2350767037116
$$
$$
-9\,{\frac { \left( n+2 \right) ^{2}}{ \left( n+3 \right) ^{2}}}+{\frac { \left( 19\,{n}^{2}+40\,n+18 \right)  \left( n+2 \right) ^{2}N}{ \left( n+3
 \right) ^{2} 
\left( n+1 \right) ^{2}}}
$$
$$
-{\frac { \left( 45+148\,n+159\,{n}^{2}+70\,{n}^{3}+11\,{n}^{4} \right) {N}^{2}}{ \left( n+3 \right) ^{2} \left( n+1
 \right) ^{2}}}+{N}^{3}
$$

$$
{\frac {9}{128}}\,{9}^{n} \left( 1+3/4\,{n}^{-1}+{\frac {53}{32}}\,{n}^{-2}+{\frac {261}{64}}\,{n}^{-3} \right) \sqrt {3}{n}^{-2}{\pi }^{-1}
$$

{\bf $(k,l)=(2,2)$, $z=1$}

$$
1, 2, 4, 10, 26, 76, 232, 764, 2578, 9076, 32264, 117448, 428936, 1589680, 5897504, 22101304, 82851218, 
$$
$$
312935236, 1182083272, 4491680504, 17067914056, 65167445872
$$

$$
128\,{\frac {n \left( n-1 \right) }{ \left( n+5 \right)  \left( n+4 \right) }}-32\,{\frac { \left( -1+4\,n+6\,{n}^{2} \right) N}{ \left( n+5 \right) 
 \left( n+4 \right) }}+
$$
$$
8\,{\frac { \left( 4+21\,n+11\,{n}^{2} \right) {N}^{2}}{ \left( n+5 \right)  \left( n+4 \right) }}-4\,{\frac { \left( -24-7\,n+{n}^{2
} \right) {N}^{3}}{ \left( n+5 \right)  \left( n+4 \right) }}-2\,{\frac { \left( 3\,n+10 \right) {N}^{4}}{n+4}}+{N}^{5}
$$

$$
1/4\,{\frac {{4}^{n}}{\pi \,n}}
$$
(for more refined asymptotics see the webpage).

{\bf $(k,l)=(2,2)$, $z=2$}

$$
1, 2, 6, 24, 120, 720, 5040, 40320, 361116, 3540600, 37207368, 411988896, 4747167568, 
$$
$$
56428884512, 687793860000, 8559142303296, 108400653865572
$$

$$
-512\,{\frac { \left( 2\,n-1 \right)  \left( n-1 \right) n}{ \left( n+5 \right)  \left( n+4 \right) ^{2}}}+16\,{\frac { \left( 231\,n+711\,{n}^{2}+164\,{n}^
{4}
+588\,{n}^{3}-2 \right) N}{ \left( n+5 \right)  \left( n+1 \right)  \left( n+4 \right) ^{2}}}
$$
$$
-16\,{\frac { \left( 1146\,{n}^{2}+545\,{n}^{3}+1003\,n+94\,
{n}^{4}+296 \right) {N}^{2}}{ \left( n+5 \right)  \left( n+1 \right)  \left( n+4 \right) ^{2}}}+4\,{\frac { \left( 2030\,{n}^{2}+85\,{n}^{4}+688\,{n}^{3}+
2501\,n+996 \right) {N}^{3}}{ \left( n+5 \right)  \left( n+1 \right)  \left( n+4 \right) ^{2}}}
$$
$$
-4\,{\frac { \left( 87\,{n}^{3}+8\,{n}^{4}+525\,n+336\,{n}^{2
}+250 \right) {N}^{4}}{ \left( n+5 \right)  \left( n+1 \right)  \left( n+4 \right) ^{2}}}+{N}^{5}
$$

$$
1/32\,{16}^{n} \left( {n}^{-1} \right) ^{7/2} \left( 1+{\frac {33}{8}}\,{n}^{-1}+{\frac {2145}{128}}\,{n}^{-2}+{\frac {81723}{1024}}\,{n}^{-3} \right) {\pi 
}^{-3/2}
$$

\end{document}